\documentclass[11pt]
{amsart}
\usepackage{hyperref}
\usepackage{pdfpages}
\hypersetup{nesting=true,debug=true,naturalnames=true}
\usepackage{graphicx,amssymb,upref}

\newtheorem{theorem}{Theorem}
\newtheorem{lemma}[theorem]{Lemma}
\newtheorem{definition}{Definition}
\title{Generalized three and four person Hat Game}
\author{Theo van Uem}  
\address{Amsterdam University of Applied Sciences, Amsterdam, The Netherlands.} 
\email{tjvanuem@gmail.com}  

\begin{document}
\hbadness=99999

\begin{abstract}
Ebert's hat problem with two colors and equal probabilities has, remarkable,  the same optimal winning probability for three and four players.
This paper studies Ebert's hat problem for three and four players, where the probabilities of the two colors may be different for each player. Our goal is to maximize the probability of winning the game and to describe winning strategies. We obtain different results for games with three and four players. We use the concept of an adequate set. The construction of adequate sets is independent of underlying probabilities and we can use this fact in the analysis of our general case.
 The computational complexity of the adequate set method is dramatically lower than by standard methods.
\end{abstract}
\maketitle
\section{Introduction}
Hat puzzles were formulated at least since Martin Gardner’s 1961 article~\cite{MG}. They have got an impulse by Todd Ebert in his Ph.D. thesis in 1998~\cite{TE}. Buhler~\cite{JB} stated: “It is remarkable that a purely recreational problem comes so close to the research frontier”. Also articles in The New York Times~\cite{SR}, Die Zeit~\cite{WB} and abcNews~\cite{JP} about this subject got broad attention. This paper studies generalized Ebert’s hat problem for three and four players. The probabilities of the two colors may be different for each player, but  known to all the players. All players guess simultaneously the color of their own hat observing only the hat colors of the other players. It is also allowed for each player to pass: no color is guessed. The team wins if at least one player guesses his or her hat color correctly and none of the players has an incorrect guess. Our goal is to maximize the probability of winning the game and to describe winning strategies.
The symmetric two color hat problem (equal probability $0.5$ for each color) with $N=2^k-1$ players is solved in~\cite{EMV}, using Hamming codes, and with $N=2^k$ players in~\cite{GC} using extended Hamming codes. 
Burke et al.~\cite{EB} try to solve the symmetric hat problem with $N=3,4,5,7$ players using genetic programming. Their conclusion: The $N$-prisoners puzzle (alternative names: Hat Problem, Hat Game) gives evolutionary computation and genetic programming a new challenge to overcome. 
Lenstra and Seroussi~\cite{HL} show that in the symmetric case of two hat colors, and for any value of $N$, playing strategies are equivalent to binary covering codes of radius one.

Krzywkowski~\cite{MK} describes applications of the hat problem and its variations, and their connections to different areas of science. 
Johnson~\cite{BJ}  ends his presentation with an open problem:
If the hat colors are not equally likely, how will the optimal strategy be affected?
We will answer this question and our method gives also interesting results in the symmetric case.
In section 2 we define our main tool: an adequate set. 
In sections 3 and 6 we obtain results for three and four person two color Hat Game, where each player $i$ may have different probabilities $(p_i,q_i)$ to get a specific colored hat.
In sections 4 and 7 we obtain results for the asymmetric three and four person two color Hat Game, where each player has the same set of probabilities $(p,q)$ to get a specific colored hat, but the probabilities are different. In sections 5 and 8 we find old and new results for the well known symmetric case $p=q=\frac{1}{2}$.
Section 9 gives a comparison between generalized three and four person Hat Game. Section 10 handles with computational complexity. Central in all our investigations are adequate sets.

\section{Adequate set}
In this section we have $N$ players and $q$ colors.
The $N$ persons in our game are distinguishable, so we can label them from 1 to
$N$. We label the $q$ colors $0,1,..,q-1.$ The probabilities of the colors are known to all players. The probability that color $i$ will be on a hat of player $j$ is
$p_{i,j}$ $(i\in\{0,1,..,q-1\},\ \ \forall j \in \{1,2,..,N\}: \sum_{i=0}^{q-1}p_{i,j}=1).$
Each possible configuration of the hats can be represented by an element of
$B=\{b_1b_2\dots b_N\vert b_i\in\left\{0,1,\dots,q-1\right\},\ i=1,2..,N\}$.
The S-code represents what the $N$ different players sees. Player $i$ sees q-ary code $b_1..b_{i-1}b_{i+1}..b_N$ with decimal value
$s_i=\sum_{k=1}^{i-1}b_k.q^{N-k-1}+\sum_{k=i+1}^Nb_k.q^{N-k}\ $, a value between 0 and $q^{N-1}-1.$ \\
Let S be the set of all S-codes:
$S=\{s_1s_2\dots s_N\vert{}s_i=\sum_{k=1}^{i-1}b_k.q^{N-k-1}+\sum_{k=i+1}^Nb_k.q^{N-k},b_i\in{}\{0,1,\dots,q-1\},\
i=1,2,\dots,N\ \}$.
Each player has to make a choice out of $q+1$ possibilities: 0='guess color 0', 
1='guess color 1', \ldots{}.,$\ q-1$ ='guess color $q-1$', $q$='pass'.
\\
We define a decision matrix $D=\left(a_{i,j}\right) \ $ where
$i\in{}\{1,2,..,N\}$(players); $j\in{}\{0,1,..,q^{N-1}-1\}$(S-code of a player);
$a_{i,j}\in{}\left\{0,1,..,q\right\}.$
The meaning of $a_{i,j}$ is: player $i$ sees S-code $j$ and takes decision $a_{i,j}$ (guess a color or pass).
We observe the total probability ($sum$) of our guesses. 
For each $b_1b_2\dots b_N$ in B   we have:\newline

IF $a_{1{,s}_1}\in{}\{q,b_1\}$  AND $a_{2,s_2}\in{}\{q,b_2\}$  AND ... AND 
$a_{N,s_N}\in{}\{q,b_N\}$  AND \\
NOT ($a_{1,s_1}=a_{2{,s}_2}=\dots=a_{N,s_N}=q$) THEN \\ $sum=sum+p_{b_1,1}.p_{b_2,2}\dots p_{b_N,N}$.
\newline

Any choice of the $a_{i,j}$ in the decision matrix determines which CASES $b_1b_2\dots b_N$ have a
positive contribution to $sum$ (we call it a GOOD CASE) and which CASES don't
contribute positive to $sum$ (we call it a BAD CASE).
\begin{definition}
 Let $A \subset B$. $A$ is adequate to $B-A$ if for each q-ary element $x$ in $B-A$ there are $q-1$ elements in A which are equal to $x$ up to one fixed q-ary position.
 \end{definition}
 
\begin{theorem}
 BAD CASES are adequate to GOOD CASES.
\end{theorem}
\begin{proof}
  Any  GOOD CASE  has at least one $a_{i,j}$ not equal to $q$. Let this
specific $a_{i,j}$ have value $b_{i_0}.$ Then our GOOD CASE generates $q-1\ $BAD CASES
by only changing  the value $b_{i_0}$ in any value of $0,1,..,q-1\ $
except $b_{i_0}$.  
\end{proof}
The notion of  an adequate set is the same idea as the concept of 
strong covering, introduced by Lenstra and Seroussi~\cite{HL}. The number of elements in an adequate set will be written as \textit{das} (dimension of adequate set).
Adequate sets are generated by an adequate set generator (ASG). See Appendix \ref{appendix:one} for an implementation in a VBA/Excel program.
Given an adequate set, we obtain a decision matrix $ D=\left(a_{i,j}\right)$ by the following procedure.\\ 

Procedure DMG (Decision Matrix Generator):\\
Begin Procedure\\
For each element in the adequate set:
\begin{itemize}
\item Determine the q-ary representation $b_1b_2\dots b_N$
 \item Calculate S-codes $s_i=\sum_{k=1}^{i-1}b_k.q^{N-k-1}+\sum_{k=i+1}^Nb_k.q^{N-k}$( $i$=1,..,N)
\item For each player $i$: fill decision matrix with  $a_{{i,s}_i}=b_i$  ($i$=1,..,N), where each cell may contain several values.
\end{itemize}
Matrix $D$ is filled with BAD COLORS. We can extract the GOOD COLORS by considering all $a_{i,j}$ with $q-1$ different BAD COLORS and then choose the only missing
color. In all situations with less than $q-1$ different BAD COLORS we pass.
When there is an $a_{i,j}$ with $q$ different BAD COLORS all colors are bad, so the first
option is to pass. But when we  choose any color, we get a situation with
$q-1$ colors.  So in case of $q$ BAD COLORS we are free to choose any color or
pass.
The code for pass is $q$, but in our decision matrices we prefer a blank, which
supports readability.
The code for `any color or pass will do' is defined  $q+1$, but we prefer a "$\star$" for readability.
\\
End Procedure.\\
 
 This procedure is implemented in the VBA/Excel program DMG (See Appendix \ref{appendix:two}).
\section{Generalized three person two color Hat Game}
{\raggedright
Three distinguishable players are randomly fitted with a white (code 0) or black (code 1) hat.
Each player $i$ has his own probabilities $p_i$ and $q_i$ to get a white
respectively a black hat, where $0<p_i<1,\  \ p_i+q_i=1\ (i=1,2,3).$ All
probabilities are known to all players.
}
Part of the strategy is that the players give themselves an identification: 1, 2 and 3.

{\raggedright
Our goal is to maximize the probability of winning the game and to describe
winning strategies.
}
Let $X$ be an adequate set and $P(X)$ is the probability generated by the adequate set.
The adequate set $X$ dominates the adequate set $Y$ if $P(X) \leq P(Y)$. We also define: $X$ dominates the adequate set $Y$  absolutely if $P(X) < P(Y)$. We use the abbreviation DOM for domination. An adequate set $A$ is non-dominated by a collection $C$ of adequate sets when $A$ dominates each element of $C$.
Adequate sets $X$ and $Y$ are isomorphic when there is a bijection from \{1,2,3\} to itself  which transforms $X$ into $Y$. The decision matrices are then also isomorphic.
{\raggedright
A player \textit{i }with $p_i<q_i$ gets an asterix: when observing such a player
we have to flip the colors: white becomes black and vice versa.
}

{\raggedright
In such a way we have without loss of generality $p_i\geq{}q_i\
\left(i=1,2,3\right).$
}

{\raggedright
The next step is to renumber the players in such a way that
$\frac{p_1}{q_1}\geq{}\frac{p_2}{q_2}\geq{}\frac{p_3}{q_3}$ , which is equivalent
to $p_1\geq{}p_2\geq{}p_3$.
}

{\raggedright
So: $1>p_1\geq{}p_2\geq{}p_3\geq{}\frac{1}{2}$
}
or, equivalently: $0<q_1\leq q_2\leq q_3\leq \frac{1}{2}.$ \\
We define decision matrices:\\
$\alpha$:
\begin{center}
     \begin{tabular}{|cccc|}
     \hline
00&01&10&11\\
\hline 
    1    &       &       & 0 \\
     1    &       &       & 0\\
      1   &       &       &  0\\
  \hline
    \end{tabular}%
  \end{center}
 $\delta$:
\begin{center}
     \begin{tabular}{|cccc|}
     \hline
00&01&10&11\\
\hline 
    0    &       &       & 1 \\
         &    0   &    1   & \\
         &    0   &   1    &  \\
  \hline
    \end{tabular}%
  \end{center}
$\epsilon:$  
\begin{center}  
     \begin{tabular}{|cccc|}
     \hline
00&01&10&11\\
\hline 
    0     & 0      & 0      & 0 \\
         &       &    $\star$   & $\star$\\
        &       &  $\star$     &$\star$  \\
  \hline
    \end{tabular}%
    \end{center}
\begin{theorem}
\label{Th1}
Given $0<q_1\leq q_2\leq q_3\leq \frac{1}{2}$ we have:\\

\begin{center}
     \begin{tabular}{|c|c|c|}
     \hline
\rule{0pt}{4ex}    CASE&$q_2=\frac{1}{2}$&$q_2<\frac{1}{2}$\\[2ex]
\hline 
   \rule{0pt}{4ex} $\frac{1}{q_1}>\frac{1}{q_2}+\frac{1}{q_3}$    & $p_1>\frac{3}{4}$ & $p_1>\frac{3}{4}$ \\
   &$\epsilon$&$\epsilon$\\
   \hline
 \rule{0pt}{4ex}$\frac{1}{q_1}<\frac{1}{q_2}+\frac{1}{q_3}$       & $\frac{3}{4}$ &  $p_1+q_1q_2q_3(\frac{1}{q_2}+\frac{1}{q_3}-\frac{1}{q_1})>p_1>\frac{3}{4}$  \\
 &$\alpha,\delta$&$\delta$\\
 \hline
\rule{0pt}{4ex}  $\frac{1}{q_1}=\frac{1}{q_2}+\frac{1}{q_3}$          &  $\frac{3}{4}$&$p_1>\frac{3}{4}$       \\
 &   $\alpha,\delta, \epsilon$  &$\delta, \epsilon$  \\ 
  \hline
    \end{tabular}%
  \end{center} 
 where in each case we give the optimal probability in the first line and the optimal decision matrices in the second line.
 \newline
 \end{theorem}
\begin{proof}
 We shall show that the following sets are absolute dominant:\\
• $\{4, 5, 6, 7\}$ when $\frac{1}{q_1}>\frac{1}{q_2}+\frac{1}{q_3}$ ,\\
• $\{3, 4\}$ and $\{0, 7\}$ when $\frac{1}{q_1}<\frac{1}{q_2}+\frac{1}{q_3}$ with $p_2=\frac{1}{2}$,\\
• $\{3, 4\}$  when $\frac{1}{q_1}<\frac{1}{q_2}+\frac{1}{q_3}$ with $p_2>\frac{1}{2}$,\\
• $\{4, 5, 6, 7\}$, $\{3, 4\}$ and $\{0, 7\}$ when $\frac{1}{q_1}=\frac{1}{q_2}+\frac{1}{q_3}$ with $p_2=\frac{1}{2}$,\\
• $\{4, 5, 6, 7\}$ and $\{3, 4\}$  when $\frac{1}{q_1}=\frac{1}{q_2}+\frac{1}{q_3}$ with $p_2>\frac{1}{2}$.\\
Let \textit{das} be the dimension of an adequate set (number of elements in the
set). Obviously, there are no adequate sets with $das<2$. When \textit{das}=2, we find 4 adequate sets (use ASG), independent of the underlying
probabilities: $ \{0,7\},\{1,6\},\{2,5\}$ and $\{3,4\}$.
We notice that $ \{1,6\},\{2,5\}$ and $\{3,4\}$ are isomorphic:
they can be obtained from any of the three by renumbering the players. E.g. interchanging players 1 and 3 in binary codes of $\left\{1,6\right\}$ gives $\left\{3,4\right\}$ and interchanging players 2 and 3 in $\left\{1,6\right\}$ gives $\left\{2,5\right\}$.  

{\raggedright
 We are looking for optimal adequate sets. An adequate set consist
of BAD CASES. We want to maximize the winning probability, so we minimize the
adequate set probability.
}

{\raggedright
The next table shows the 4 adequate sets and probabilities:
}

\begin{center}
     \begin{tabular}{|c|c|}
     \hline
\{0,7\} &$p_1p_2p_3+q_1q_2q_3=A$ \\
  \hline
\{1,6\} &$p_1p_2q_3+q_1q_2p_3=B$ \\
  \hline
 \{2,5\} &$p_1q_2p_3+q_1p_2q_3=C$ \\
  \hline
 \{3,4\} &$p_1q_2q_3+q_1p_2p_3=D$ \\
  \hline
    \end{tabular}%
\end{center}

{\raggedright
We have:
}

\begin{align*}
A-B&=q_1q_2q_3(\frac{p_1}{q_1}\frac{p_2}{q_2}-1)(\frac{p_3}{q_3}-1)\\
A-C&=q_1q_2q_3(\frac{p_1}{q_1}\frac{p_3}{q_3}-1)(\frac{p_2}{q_2}-1)\\
A-D&=q_1q_2q_3(\frac{p_3}{q_3}\frac{p_2}{q_2}-1)(\frac{p_1}{q_1}-1)\\
B-C&=q_1q_2q_3(\frac{p_2}{q_2}-\frac{p_3}{q_3})(\frac{p_1}{q_1}-1); \\  
B-D&=q_1q_2q_3(\frac{p_1}{q_1}-\frac{p_3}{q_3})(\frac{p_2}{q_2}-1) \\  
C-D&=q_1q_2q_3(\frac{p_1}{q_1}-\frac{p_2}{q_2})(\frac{p_3}{q_3}-1)
\end{align*}

{\raggedright
So we have: $A\geq{}B\geq{}C\geq{}D$: the adequate set \{3,4\} dominates all
other adequate sets when \textit{das}=2.
}

{\raggedright
When \textit{das}$=3$, we get (using the adequate set generator) 24 adequate sets (see Appendix \ref{appdas3}),
all absolutely dominated by $ \{0,7\},\{1,6\},\{2,5\}$ or $\{3,4\}$.
}\\
When \textit{das}=4 we get the situation in Appendix \ref{appendix:seven}.
\newline
When $das>4$ then always one or more of the sets $ \{0,7\},\{1,6\},\{2,5\}$, $\{3,4\}$ is included
 (we don't need the ASG).\\

Adequate set $\{0,7\}$ has value A and decision matrix $\alpha$. Adequate set $\{3,4\}$ has value D and decision matrix $\delta$.  Adequate set $\{4,5,6,7\}$ has value $p_1$ and decision matrix $\epsilon$.
 Using DMG (Appendix \ref{appendix:two}) we find $\alpha,  \delta$ and $\epsilon$.
 \\ [2ex]
 We first consider the battle between \{3,4\} and \{4,5,6,7\}.\\
$\{3,4\}$ is the winner when $p_1q_2q_3+q_1p_2p_3 <
q_1$, so:
 $\frac{1}{q_1}<\frac{1}{q_2}+\frac{1}{q_3}$\\
 \{4,5,6,7\} is the winner when $\frac{1}{q_1}>\frac{1}{q_2}+\frac{1}{q_3}$\\
 CASE 1 \\
$\frac{1}{q_1}>\frac{1}{q_2}+\frac{1}{q_3}$ :\\
We have $\frac{1}{q_1}> \frac{1}{q_2}+\frac{1}{q_3}$ , so $q_1 < q_2$; we get: \\  $ 1>p_1>p_2\geq p_3\geq \frac{1}{2} \geq q_3 \geq q_2 >q_1>0$.\\
Consider Appendix \ref{appendix:seven}. To obtain absolute dominance we have to examine the adequate sets $\{0, 3, 5, 6\}$ and $\{1, 2, 4, 7\}$.
For the set $\{0, 3, 5, 6\}$, the probability is $q_1 + (p_1-q_1)(p_2p_3 + q_2q_3)
> q_1$.\\
Similarly, for the set $\{1, 2, 4, 7\}$, its probability is also greater than $q_1$.

So $ \{4,5,6,7\}$ is absolute dominant with winning probability $1-q_1=p_1$ and decision matrix $\epsilon$, where $\frac{1}{q_1}>\frac{1}{q_2}+\frac{1}{q_3}\geq 4$, so $ p_1>\frac{3}{4}$ 
\\ [2ex]
CASE 2 \\
$\frac{1}{q_1}<\frac{1}{q_2}+\frac{1}{q_3}$ \\
There are 4 potential dominant sets: $\{0,7\},\{1,6\},\{2,5\},\{3,4\}$.
The last three are isomorphic and because of $A\leq B\leq C\leq D$ we analyze the battle between $A$ and $D$.
\\ [2ex]
CASE 2.1. $\alpha $ and      $\delta$ are both optimal: $A=D$, so: 
\begin{equation*}
(\frac{p_3}{q_3}\frac{p_2}{q_2}-1)(\frac{p_1}{q_1}-1)=0 \Leftrightarrow 
\end{equation*}
\begin{equation*}
((p_2=q_2=\frac{1}{2}) \land (p_3=q_3= \frac{1}{2})) \lor (p_1=q_1=\frac{1}{2}) \Leftrightarrow
\end{equation*}
\begin{equation*}
(\frac{1}{2}=p_3=p_2\leq p_1<1)\lor (p_1=p_2=p_3=\frac{1}{2}) \Leftrightarrow
\end{equation*}
\begin{equation*}
(\frac{1}{2}=p_3=p_2\leq p_1<1) \Leftrightarrow 
\end{equation*}
\begin{equation*}
p_2=\frac{1}{2} 
\end{equation*}
and  optimal probability is $1-(p_1q_2q_3+q_1p_2p_3)=\frac{3}{4}$.\\ [2ex]
CASE 2.2. Only $\delta$ optimal: $A<D$, so:
$p_2>\frac{1}{2}$.\\ We note:
 $ \{1,6\},\{2,5\}$ and $\{3,4\}$ are isomorphic, but this doesn't imply  $B=C=D$: the probabilities are not always invariant by a bijection of the players.
We have $A<D$ and $B\leq C\leq D$. When $B$ or $C$ are equal to $D$ then D is absolute dominant (up to isomorphic) otherwise $B\leq  C< D$ and D is absolute dominant.
Optimal probability: $1-(p_1q_2q_3+q_1p_2p_3)=p_1+q_1q_2q_3(\frac{1}{q_2}+\frac{1}{q_3}-\frac{1}{q_1})>p_1$.  And $\frac{1}{q_1}=\frac{1}{q_2}+\frac{1}{q_3}<4$, so $ p_1>\frac{3}{4}$  \\ [2ex]
CASE 3 \\
$\frac{1}{q_1}=\frac{1}{q_2}+\frac{1}{q_3}$ \\ 
Optimal probability: $p_1+q_1q_2q_3(\frac{1}{q_2}+\frac{1}{q_3}-\frac{1}{q_1})=p_1$\\
CASE 3.1 $p_2=\frac{1}{2}$\\
Optimal probability: $p_1=\frac{3}{4}$\\ 
Optimal decision matrices:\\
$\frac{1}{q_1}\geq \frac{1}{q_2}+\frac{1}{q_3}$ gives $\epsilon$ and $\frac{1}{q_1} \leq \frac{1}{q_2}+\frac{1}{q_3}$ gives $\alpha $ and $
\delta$.\\
CASE 3.2 $p_2>\frac{1}{2}$\\
Optimal probability: 
$p_1$ and $\frac{1}{q_1}=\frac{1}{q_2}+\frac{1}{q_3}<4$, so $ p_1>\frac{3}{4}$ \\
Optimal decision matrices:\\
$\frac{1}{q_1}\geq \frac{1}{q_2}+\frac{1}{q_3}$ gives $\epsilon$ and $\frac{1}{q_1} \leq \frac{1}{q_2}+\frac{1}{q_3}$ gives  $
\delta$.
\end{proof} 

Note 1: Instead of $\frac{1}{q}$ we can also use $\frac{p}{q}$. E.g. : we get $\frac{p_1}{q_1}=\frac{p_2}{q_2}+\frac{p_3}{q_3}+1$ instead of $\frac{1}{q_1}=\frac{1}{q_2}+\frac{1}{q_3}$.\\
\\
Note 2: We observe that the well-known strategy $\alpha$ is only dominant when \\      
$(\frac{1}{q_1}\leq \frac{1}{q_2}+\frac{1}{q_3}) \land (q_2=\frac{1}{2}) \Leftrightarrow (p_2=p_3=\frac{1}{2}) \land (\frac{1}{2}\leq p_1\leq \frac{3}{4}).$ 
\section{Asymmetric three person two color Hat Game }
In this section we study three person two color asymmetric Hat Game. For each player let $p$ be the
probability to get a white hat and $q$ be the probability to get a black hat.
Without loss of generality we may assume (asymmetric case): $\frac{1}{2}<p<1$.

\begin{theorem}
In asymmetric three person two color hat game we have maximal probability $1-pq$ of winning the game, with decision matrix:
\begin{center}
     \begin{tabular}{|cccc|}
     \hline
00&01&10&11\\
\hline 
       0  &       &       & 1 \\
        &   0    &   1    & \\
         &  0     &  1     &  \\
  \hline
    \end{tabular}%
\end{center} 

\end{theorem}
\begin{proof}
Use the result of CASE $(\frac{1}{q_1}<\frac{1}{q_2}+\frac{1}{q_3}) \land (p_2>\frac{1}{2})$ in Theorem \ref{Th1}.
Optimal probability is $1-(p_1q_2q_3+q_1p_2p_3)=1-pq$.    
 \end{proof}
    
\section{Symmetric two color three person Hat Game}
In this section we focus on the symmetric Hat Game with two colors and three players. Each player has a white hat with probability $\frac{1}{2}$ and a black hat with probability $\frac{1}{2}$. 
\begin{theorem}
For symmetric three person two color Hat Game the
 maximal probability is $\frac{3}{4}$ and the optimal decision matrices are:
\begin{center}
     \begin{tabular}{|cccc|}
     \hline
00&01&10&11\\
\hline 
   0      &       &       & 1 \\
        &   0   &   1    & \\
          & 0     &  1    &  \\
  \hline
    \end{tabular}%
\end{center} 

and:
\begin{center}
     \begin{tabular}{|cccc|}
     \hline
00&01&10&11\\
\hline 
    1     &       &       & 0 \\
    1     &       &       & 0 \\
    1     &       &       & 0 \\
  \hline
    \end{tabular}%
\end{center} 
\end{theorem}
\begin{proof}
Use result in CASE  $(\frac{1}{q_1}<\frac{1}{q_2}+\frac{1}{q_3}) \land (p_2=\frac{1}{2})$ in Theorem \ref{Th1}.
\end{proof}

\section{Generalized four person two color Hat Game}
{\raggedright
Four distinguishable players are randomly fitted with a white or black hat. The code for a white hat is 0 and for a black hat is 1.
Each player $i$ has his own probabilities $p_i$ and $q_i$ to get a white
respectively a black hat, where $0<p_i<1,\  \ p_i+q_i=1\ (i=1,2,3,4).$ All
probabilities are known to all players.
}
Part of the strategy is that the players give themselves an identification: 1, 2, 3  and 4.

{\raggedright
Our goal is to maximize the probability of winning the game and to describe
winning strategies.
}

{\raggedright
A player \textit{i }with $p_i<q_i$ gets an asterix: when observing such a player
we have to flip the colors: white becomes black and vice versa.
}

{\raggedright
In such a way we have without loss of generality $p_i\geq{}q_i\
\left(i=1,2,3,4\right).$
}

{\raggedright
The next step is to renumber the players in such a way that
$\frac{p_1}{q_1}\geq{}\frac{p_2}{q_2}\geq{}\frac{p_3}{q_3}\geq \frac{p_4}{q_4}$ , which is equivalent
to $p_1\geq{}p_2\geq{}p_3\geq p_4$.
}

{\raggedright
So: $1>p_1\geq{}p_2\geq{}p_3\geq{}p_4\geq \frac{1}{2}$
}
or, equivalently: $0<q_1\leq q_2\leq q_3\leq q_4\leq \frac{1}{2}.$ \\
Using ASG  with \textit{das}$<4$ we get no adequate sets.  
\begin{lemma}
 $\{6, 7, 8, 9\}$ dominates all adequate sets with \textit{das}=4
\end{lemma}
\begin{proof}
When \textit{das}=4 we get, using ASG,  40 adequate sets, see appendix \ref{6789}. In the same Appendix we see that all sets are dominated by another set, except $\{6,7,8,9\}$.
The proof of each dominance relation proceeds along the same lines. We give one example: nr 2. dominates nr 1.\\
$\{0,2,13,15\}$ DOM $\{0,1,14,15\}$\\
1, 14 $\geq$ 2,13 \quad 
0001 1110 $\geq$ 0010 1101 \\
$(\frac{p_1}{q_1}.\frac{p_2}{q_2}-1)(\frac{p_1}{q_1}-\frac{p_2}{q_2})\geq 0$.
\end{proof}

\begin{lemma}
When \textit{das}=5  we get 8 non-dominated adequate sets:
\begin{center}
\begin{tabular}{ |c|c|c|c|c|c| } 
 \hline
 {$S_0$} & 0& 11 &13  & 14 &15 \\ 
 {$S_1$} & 1 & 10 &12  &14  &15 \\
{$S_2$}  & 2 & 9 &12  &13  &15 \\
 {$S_3$}   &3  & 8 & 12 & 13 &14 \\
  {$S_4$}   & 4 &9 & 10 & 11 &15 \\
   {$S_5$}   & 5 & 8 & 10 &11  &14 \\
    {$S_6$}   & 6 & 8 &9  &11  &13 \\
     {$S_7$}   &7  & 8 & 9 & 10 &12 \\
       
 \hline
\end{tabular}
\end{center}
where $S_x$ is the shortcut for the adequate set starting with an $x$.
\end{lemma}
\begin{proof}
When \textit{das}=5 we get, using ASG, 560 adequate sets. 
We use the inclusion principle: we look for subsets of these 560 sets in the set of 40 adequate sets with \textit{das}=4. This procedure is realized in the program ASG45 (see appendix \ref{appendix:four}). In this way we eliminate 480 dominated adequate sets. The output of ASG45 (80 non-dominated sets) is shown in appendix \ref{appendix:five}, where we also show that $S_x (x=0,1,..,7)$ are the only non-dominated sets when \textit{das}=5.
\end{proof}
\begin{lemma}
When \textit{das} $>5$  we find one non-dominated adequate  set:\\ $\{8, 9, 10, 11, 12, 13, 14, 15\}$
\end{lemma}
\begin{proof}
When \textit{das}$>$5 then only dominated sets are found, with one exception: \textit{das}=8. Use the inclusion principle: a subset is found in the \textit{das} = 5 adequate sets; this procedure can be automated in programs ASG56, ASG57,...ASG516, analogue to program ASG45.
Using ASG58 we get  10 (out of 10310) non-dominated sets :
\begin{center}
\begin{tabular}{ |c|c|c|c|c|c|c|c| } 
 \hline
0&	1&	2	&3	&4&	5	&6	&7\\
0&	1&	2&	3	&8&9	&10&	11\\
0&	1&	4	&5&8	&9&	12&	13\\
0	&2&	4	&6&	8&	10&	12&	14\\
0&	3&	5&	6&	9&	10&	12&	15\\
1&	2&	4&	7	&8&	11&	13&	14\\
1&	3&	5&	7&	9&	11&	13&	15\\
2&	3&	6&	7&	10&	11&	14&	15\\
4&	5&	6&	7&	12&	13&	14&	15\\
8&	9&	10&	11&	12&	13&	14&	15\\
\hline
\end{tabular}
\end{center}
They are all dominated (use from top to bottom the four different positions principle) by the last one.
For example first and second element:
$\{0,1,2,3,8,9,10,11\}$ DOM $\{0,1,2,3,4,5,6,7\}$ \\
4,5,6,7 $\geq $ 8,9,10,11\\
0100 0101 0110 0111 $\geq$ 1000 1001 1010 1011 \\
(01-10)00+(01-10)01+(01-10)10+(01-10)11 $\geq 0$\\
which results in: $\frac{p_1}{q_1}-\frac{p_2}{q_2}\geq 0$.
\\ We abbreviate the last one by its first element: $S_8$.
\end{proof}
The optimal set when \textit{das}=4: $\{6,7,8,9\}$ is noted by its last element: $S_9$.\\ 
Our goal is to prove: $\{S_7,S_8,S_9\}$ dominates all adequate sets.\\
We make use of the following Lemmas:
\begin{lemma}
$\{S_7,S_9\}$ dominates $S_6$
\end{lemma}
\begin{proof}
$S_7$= $\{7,8,9,10,12\}$ dominates $S_6$=$\{6,8,9,11,13\}$  when: \\
7 10 12 $\leq$ 6 11 13\\
0111 1010 1100 $\leq$ 0110 1011 1101\\
$(-p_1q_2q_3+q_1p_2q_3+q_1q_2p_3)(p_4-q_4)\leq 0$\\
$(-\frac{p_1}{q_1}+\frac{p_2}{q_2}+\frac{p_3}{q_3}) (\frac{p_4}{q_4}-1)\leq 0$\\

$S_9$=$\{6,7,8,9\}$  dominates $S_6$= $\{6,8,9,11,13\}$ when: \\
7  $\leq$  11 13\\
0111 $\leq$ 1011 1101\\
$(p_1q_2q_3-q_1p_2q_3-q_1q_2p_3)q_4\leq 0$\\
$\frac{p_1}{q_1}-\frac{p_2}{q_2}-\frac{p_3}{q_3} \leq 0$\\
\end{proof}

\begin{lemma}
$\{S_7,S_9\}$ dominates $S_5$ 
\end{lemma}
\begin{proof}
$S_7$=$\{7,8,9,10,12\}$  dominates $S_5$= $\{5,8,10,11,14\}$ when: \\
7 9 12 $\leq$ 5 11 14\\
0111 1001 1100 $\leq$ 0101 1011 1110\\

$(-\frac{p_1}{q_1}+\frac{p_2}{q_2}+\frac{p_4}{q_4}) (\frac{p_3}{q_3}-1)\leq 0$\\

$S_9$= $\{6,7,8,9\}$ dominates $S_5$= $\{5,8,10,11,14\}$ when: \\
 6 7 9 $\leq$ 5 10  11 14\\
0110 0111 1001 $\leq$ 0101 1010 1011 1110\\
$(\frac{p_1}{q_1}-\frac{p_2}{q_2})(\frac{p_3}{q_3}-\frac{p_4}{q_4}) \geq \frac{p_1}{q_1}-\frac{p_2}{q_2}-\frac{p_4}{q_4}$\\
\end{proof}

\begin{lemma}
$\{S_7,S_9\}$ dominates $S_4$
\end{lemma}
\begin{proof}
$S_7$=$\{7,8,9,10,12\}$ dominates $S_4$=$\{4,9,10,11,15\}$  when: \\
7 8 12 $\leq$ 4 11 15\\
0111 1000 1100 $\leq$ 0100 1011 1111\\
$(\frac{p_1}{q_1}-\frac{p_2}{q_2}-1) (\frac{p_3}{q_3}\frac{p_4}{q_4}-1)\geq 0$\\

$S_9$=$\{6,7,8,9\}$  dominates $S_4$= $\{5,9,10,11,15\}$ when: \\
 6 7 8 $\leq$ 4 10  11 15\\
 0110 0111 1000 $\leq$ 0100 1010 1011 1111\\
$(\frac{p_1}{q_1}-\frac{p_2}{q_2})(\frac{p_3}{q_3}-1)\frac{p_4}{q_4} \geq \frac{p_1}{q_1}-\frac{p_2}{q_2}-1$\\
\end{proof}

\begin{lemma}
$\{S_7,S_9\}$ dominates $S_3$ 
\end{lemma}
\begin{proof}
$S_7$=$\{7,8,9,10,12\}$ dominates $S_3$= $\{3,8,12,13,14\}$ when: \\
7 9 10   $\leq$ 3 13 14 \\
0111 1001 1010 $\leq$ 0011 1101 1110\\
$(\frac{p_1}{q_1}-\frac{p_3}{q_3}-\frac{p_4}{q_4})(\frac{p_2}{q_2}-1)\geq 0$\\

$S_9$= $\{6,7,8,9\}$ dominates $S_3$=$\{3,8,12,13,14\}$  when: \\
 6 7 9 $\leq$ 3 12 13 14\\
 0110 0111 1001 $\leq$ 0011 1100 1101 1110\\
$(\frac{p_1}{q_1}-\frac{p_3}{q_3})(1+\frac{q_4}{p_4}-\frac{p_2}{q_2}\frac{q_4}{p_4}) \leq 1$\\
When $\frac{p_1}{q_1}-\frac{p_3}{q_3}-\frac{p_4}{q_4}\geq 0$ then $S_7$ dominates $S_3$. So we consider $\frac{p_1}{q_1}-\frac{p_3}{q_3}-\frac{p_4}{q_4}\leq 0$ and get in the $S_9$ dominates $S_3$ case: \\
$(\frac{p_1}{q_1}-\frac{p_3}{q_3})(1+\frac{q_4}{p_4}-\frac{p_2}{q_2}\frac{q_4}{p_4}) \leq \frac{p_4}{q_4}(1+\frac{q_4}{p_4}-\frac{p_2}{q_2}\frac{q_4}{p_4})=\frac{p_4}{q_4}-\frac{p_2}{q_2}+1 \leq 1$
\end{proof}

\begin{lemma}
$\{S_7,S_9\}$ dominates $S_2$
\end{lemma}
\begin{proof}
After some calculations we get: \\
$S_7$ dominates $S_2$ when 
$(\frac{p_1}{q_1}-\frac{p_3}{q_3}-1)(\frac{p_2}{q_2}\frac{p_4}{q_4}-1) \geq 0$\\

$S_9$ dominates $S_2$ when
$(\frac{p_1}{q_1}-\frac{p_3}{q_3})(1+\frac{p_4}{q_4}-\frac{p_2}{q_2}\frac{p_4}{q_4}) \leq 1$\\
When $\frac{p_1}{q_1}-\frac{p_3}{q_3}\geq 1$  then $S_7$ dominates $S_2$. So we consider $\frac{p_1}{q_1}-\frac{p_3}{q_3}\leq 1$ and get in the $S_9$ dominates $S_2$ case: \\
$(\frac{p_1}{q_1}-\frac{p_3}{q_3})(1+\frac{p_4}{q_4}-\frac{p_2}{q_2}\frac{p_4}{q_4}) \leq 1+\frac{p_4}{q_4}(1-\frac{p_2}{q_2}) \leq 1$
\end{proof}

\begin{lemma}
$\{S_7,S_9\}$ dominates $S_1$ 
\end{lemma}
\begin{proof}
After some calculations we get: \\
$S_7$ dominates $S_1$ when 
$(\frac{p_1}{q_1}-\frac{p_4}{q_4}-1)(\frac{p_2}{q_2}\frac{p_3}{q_3}-1) \geq 0$\\

$S_9$ dominates $S_1$ when
$(\frac{p_1}{q_1}-\frac{p_4}{q_4})(1-\frac{p_2}{q_2}\frac{p_3}{q_3})+(\frac{p_1}{q_1}-\frac{p_2}{q_2})\frac{p_4}{q_4}+(\frac{p_2}{q_2}-\frac{p_4}{q_4})\frac{p_3}{q_3} \geq 1$
\\
When $\frac{p_1}{q_1}-\frac{p_4}{q_4}\geq 1$  then $S_7$ dominates $S_1$. So we consider $\frac{p_1}{q_1}-\frac{p_4}{q_4}\leq 1$ and get in the $S_9$ dominates $S_1$ case: \\
$(\frac{p_1}{q_1}-\frac{p_4}{q_4})(1-\frac{p_2}{q_2}\frac{p_3}{q_3})+(\frac{p_1}{q_1}-\frac{p_2}{q_2})\frac{p_4}{q_4}+(\frac{p_2}{q_2}-\frac{p_4}{q_4})\frac{p_3}{q_3} \leq (\frac{p_1}{q_1}-\frac{p_4}{q_4})(1-\frac{p_2}{q_2}\frac{p_3}{q_3})+(\frac{p_1}{q_1}-\frac{p_2}{q_2})\frac{p_3}{q_3}+(\frac{p_2}{q_2}-\frac{p_4}{q_4})\frac{p_3}{q_3}=(\frac{p_1}{q_1}-\frac{p_4}{q_4})(1-\frac{p_2}{q_2}\frac{p_3}{q_3}+\frac{p_3}{q_3})\leq 1+(1-\frac{p_2}{q_2})\frac{p_3}{q_3} \leq 1 $
\end{proof}
In the next  lemma we use$\{S_6,S_9\}$ instead of $\{S_7,S_9\}$, which gives no problems for $\{S_7,S_9\}$ dominates $S_6$.
\begin{lemma}
$\{S_6,S_9\}$ dominates $S_0$ 
\end{lemma}
\begin{proof}
After some calculations we get: \\
$S_6$ dominates $S_0$ when 
$(\frac{p_4}{q_4}(\frac{p_1}{q_1}-1)-1)(\frac{p_2}{q_2}\frac{p_3}{q_3}-1) \geq 0$\\

$S_9$ dominates $S_0$ when
$(\frac{p_1}{q_1}-1)(1-\frac{p_2}{q_2}\frac{p_3}{q_3})\frac{p_4}{q_4}+(\frac{p_1}{q_1}-\frac{p_2}{q_2})+(\frac{p_2}{q_2}-1)\frac{p_3}{q_3} \leq 1$
\\
When $\frac{p_4}{q_4}(\frac{p_1}{q_1}-1)\geq 1$ then $S_6$ dominates $S_0$. So we consider $\frac{p_4}{q_4}(\frac{p_1}{q_1}-1)\leq 1$ and get in the $S_9$ dominates $S_0$ case: \\
$(\frac{p_1}{q_1}-1)(1-\frac{p_2}{q_2}\frac{p_3}{q_3})\frac{p_4}{q_4}+(\frac{p_1}{q_1}-\frac{p_2}{q_2})+(\frac{p_2}{q_2}-1)\frac{p_3}{q_3}\leq (\frac{p_1}{q_1}-1)(1-\frac{p_2}{q_2}\frac{p_3}{q_3})\frac{p_4}{q_4}+(\frac{p_1}{q_1}-\frac{p_2}{q_2})\frac{p_3}{q_3}+(\frac{p_2}{q_2}-1)\frac{p_3}{q_3}=(\frac{p_1}{q_1}-1)\frac{p_4}{q_4}+(\frac{p_1}{q_1}-1)\frac{p_3}{q_3}(1-\frac{p_2}{q_2}\frac{p_4}{q_4}) \leq 1$ 
\end{proof}
All these Lemmas leads us to
\begin{theorem}
When $\frac{1}{q_1} \leq \frac{1}{q_2}+\frac{1}{q_3}-1$ then maximal winning probability is $p_1+q_1q_2q_3(\frac{1}{q_2}+ \frac{1}{q_3}-\frac{1}{q_1}) \geq p_1+q_1q_2q_3$ with optimal decision matrix:
\begin{center}
     \begin{tabular}{|cccccccc|}
     \hline
   {000} & 001     & 010    & 011    & 100   & 101   & 110   & {111} \\
\hline
          0&0      &       &      &      &       &   1   &    1 \\
          &      &   0    &     0 &    1  &  1     &      &     \\
          &       &  0     &      0 &    1   &  1     &     &     \\
         &      &       &      $\star$ &   $\star$    &       &      &  \\
 \hline
    \end{tabular}%
 \end{center}
 ($\star$ means: any color or pass; stars are independent)
 
When $ \frac{1}{q_2}+\frac{1}{q_3}-1 \leq \frac{1}{q_1} \leq \frac{1}{q_2}+\frac{1}{q_3}+\frac{1}{q_4}-1 $ then maximal winning probability is $p_1+q_1q_2q_3q_4(\frac{1}{q_2}+ \frac{1}{q_3}+\frac{1}{q_4}-\frac{1}{q_1}-1)$, a value between $p_1$ and $p_1+q_1q_2q_3$, with optimal decision matrix:
\begin{center}
     \begin{tabular}{|cccccccc|}
     \hline
   {000} & 001     & 010    & 011    & 100   & 101   & 110   & {111} \\
\hline
     0     &  0    & 0      &      &   0   &       &      &    1 \\
          &      &       &     0 & $\star$     &  1     &     1 &     \\
          &      &       &     0 & $\star$     &  1     &     1 &     \\
                    &      &       &     0 & $\star$     &  1     &     1 &     \\
 \hline
    \end{tabular}%
 \end{center}
 
When $  \frac{1}{q_1} \geq \frac{1}{q_2}+\frac{1}{q_3}+\frac{1}{q_4}-1 $ then maximal winning probability is $p_1$ with optimal decision matrix:
\begin{center}
     \begin{tabular}{|cccccccc|}
     \hline
   {000} & 001     & 010    & 011    & 100   & 101   & 110   & {111} \\
\hline
      0    & 0	     &      0 &    0  &  0    &   0    &     0 &   0  \\
          &      &       &      &  $\star$    &  $\star$     &    $\star$  & $\star$    \\
           &      &       &      &  $\star$    &  $\star$     &    $\star$  & $\star$    \\
                     &      &       &      &  $\star$    &  $\star$     &    $\star$  & $\star$    \\
 \hline
    \end{tabular}%
 \end{center}

\end{theorem}
\begin{proof}
The preceding Lemmas invite us to make a comparison between $S_7$, $S_8$ and $S_9$.
After some calculations we get:\\
\\
$S_9$ is winner when ($S_9$ dominates $S_8$) and ($S_9$ dominates $S_7$):\\ $(\frac{1}{q_1}\leq \frac{1}{q_2}+\frac{1}{q_3}) \land (\frac{1}{q_1}\leq \frac{1}{q_2}+\frac{1}{q_3}-1)=(\frac{1}{q_1}\leq \frac{1}{q_2}+\frac{1}{q_3}-1)$\\
 Maximal probability: $1-(p_1q_2q_3+q_1p_2p_3)=
 p_1+q_1q_2q_3(\frac{1}{q_2}+ \frac{1}{q_3}-\frac{1}{q_1}) \geq p_1+q_1q_2q_3$\\
 \\
$S_7$ is winner when ($S_7$ dominates $S_9$) and ($S_7$ dominates $S_8$):\\ $(\frac{1}{q_1}\geq \frac{1}{q_2}+\frac{1}{q_3}-1) \land (\frac{1}{q_1}\leq \frac{1}{q_2}+\frac{1}{q_3}+\frac{1}{q_4}-1)=(\frac{1}{q_2}+\frac{1}{q_3}-1\leq \frac{1}{q_1}\leq \frac{1}{q_2}+\frac{1}{q_3}+\frac{1}{q_4}-1)$\\
Maximal probability:  $1-{[}p_1q_2q_3q_4+q_1p_2p_3+q_1p_4(p_2q_3+q_2p_3){]}=p_1+q_1q_2q_3q_4(\frac{1}{q_2}+ \frac{1}{q_3}+\frac{1}{q_4}-\frac{1}{q_1}-1)$, a value between $p_1$ and $p_1+q_1q_2q_3$.\\
\\
$S_8$ is winner when ($S_8$ dominates $S_9$) and ($S_8$ dominates $S_7$):\\ $(\frac{1}{q_1}\geq \frac{1}{q_2}+\frac{1}{q_3}) \land (\frac{1}{q_1}\geq \frac{1}{q_2}+\frac{1}{q_3}+\frac{1}{q_4}-1)=(\frac{1}{q_1}\geq \frac{1}{q_2}+\frac{1}{q_3}+\frac{1}{q_4}-1)$\\
Maximal probability: $p_1$.

\end{proof}
Note 1.\\
 Conditions can also be formulated in the form $\frac{p}{q}$, e.g. $ \frac{1}{q_2}+\frac{1}{q_3}-1 \leq \frac{1}{q_1} \leq \frac{1}{q_2}+\frac{1}{q_3}+\frac{1}{q_4}-1 $ becomes $ \frac{p_2}{q_2}+\frac{p_3}{q_3} \leq \frac{p_1}{q_1} \leq \frac{p_2}{q_2}+\frac{p_3}{q_3}+\frac{p_4}{q_4}+1 $.\\

Note 2.\\
The domination used in this section is not absolutely, so there may be more optimal decision matrices with the same maximal probability. In sections 7 and 8 we shall use absolute domination and get all non isomorphic optimal decision matrices.
 
\section{Asymmetric four person Two Color Hat Game}
\label{threetwo}
\begin{theorem}
In asymmetric four person (two color) hat game we have maximal probability $1-pq$ of winning the game, with two optimal decision matrices:
\begin{center}
     \begin{tabular}{|cccccccc|}
     \hline
   {000} & 001     & 010    & 011    & 100   & 101   & 110   & {111} \\
\hline
          & 1     &       & 1     & 0     &       & 0     &     \\
          & 1     &       & 1     & 0     &       & 0     &     \\
          & $\star$      &       &       &       &       &$\star$     &     \\
    0     & 0     &       &       &       &       & 1     & 1 \\
 \hline
    \end{tabular}%
 \end{center}
 and:
 \begin{center}
     \begin{tabular}{|cccccccc|}
     \hline
   {000} & 001     & 010    & 011    & 100   & 101   & 110   & {111} \\
\hline
          & 1     & 0      &     &      &    0   & 1     &     \\
          & 1     & 0      &      &      & 0      & 1     &     \\
          & 1      &  0     &       &   0    &       &     &  1   \\
    0     &      &       &   1    &       & 1      & 0     & \\
 \hline
    \end{tabular}%
 \end{center}
 {\raggedright

}
\end{theorem}
\begin{proof}
Partial results can be obtained as a special case of the preceding section, but we want all (non isomorphic) optimal decision matrices and therefore we use absolute domination.\\
We use  ASG with parameters $n$=4, $p$=0.9, 
\textit{das}=4 and get 40 adequate sets.   Minimum sum is 0.09 and 24 adequate sets are optimal. Appendix \ref{appendix:three} shows a sorted list of all 40 adequate sets.
{\raggedright
By definition, the construction of an adequate set is independent of $p$.
}
{\raggedright
In Appendix \ref{appendix:three} we use patterns. E.g. pattern 01210 correspondents with probability $0*q^4+1*pq^3+2*p^2q^2+1*p^3q+0*p^4$. We get:
 \begin{center}
    \begin{tabular}{|l|l|}
    \hline
    pattern& probability \\
    \hline
    01210 (sum: 0.09) & $pq^3+2p^2q^2+p^3q$ \\
    10120 (sum: 0.154) & $q^4+p^2q^2+2p^3q$ \\
    02101 (sum: 0.666) & $2pq^3+p^2q^2+p^4$ \\
    11011 (sum: 0.73) & $q^4+pq^3+p^3q+p^4$ \\
    \hline
    \end{tabular}%
 \end{center}
It is not difficult to prove that 01210 absolutely dominates all other patterns when
$p$$>$$q$: the probability of 01210 is less than all the other probabilities.
}

{\raggedright
Let $\Psi{}\left(N,p\right)$ be the maximum probability of correct guessing in
our asymmetric hat game with $N$ players..
}
\[
\Psi{}\left(4,p\right)=1-\left(\ pq^3+2{\ p}^2q^2+{\ \
p}^3q\right)=1-pq=1-p+{\ p}^2
\]
{\raggedright
We remark that $\Psi{}\left(4,p\right)=\Psi{}\left(3,p\right)$ .
}\\
{\raggedright
All  24 adequate sets with the 01210 pattern generates optimal decision matrices.
Procedure DMG  gives as  result with adequate
set \{1,3,12,14\}:
 \begin{center}
     \begin{tabular}{|cccccccc|}
     \hline
   {000} & 001     & 010    & 011    & 100   & 101   & 110   & {111} \\
\hline
          & 1     &       & 1     & 0     &       & 0     &     \\
          & 1     &       & 1     & 0     &       & 0     &     \\
          & $\star$      &       &       &       &       &$\star$     &     \\
    0     & 0     &       &       &       &       & 1     & 1 \\
 \hline
    \end{tabular}%
 \end{center}
 Where $\star$ means: any color or pass will do. This happens when player 3 sees 001 or
110, which corresponds to situations 0001, 0011, 1100, 1110. In all these
situations player 1 guesses wrong, so the guess of player 3 is irrelevant.
}\\
We concentrate on the 24 optimal decision matrices and observe 12 matrices where one player can always pass and 12 matrices of a different structure. 
Appendix \ref{appendix:three} shows the two groups of 12 elements, the position of the player who can always PASS and the CYCLE to obtain isomorphic relation with the first element of each group:
 the 24 optimal adequate sets can be divided in two groups of each 12 isomorphic elements. So the first 12 rows are isomorphic to the adequate set \{1,3,12,14\} and the next 12 rows   
are isomorphic to the adequate set \{1,6,10,13\} with decision matrix:
\begin{center}
     \begin{tabular}{|cccccccc|}
     \hline
   {000} & 001     & 010    & 011    & 100   & 101   & 110   & {111} \\
\hline
          & 1     & 0      &     &      &    0   & 1     &     \\
          & 1     & 0      &      &      & 0      & 1     &     \\
          & 1      &  0     &       &   0    &       &     &  1   \\
    0     &      &       &   1    &       & 1      & 0     & \\
 \hline
    \end{tabular}%
 \end{center}
Two adequate sets are equivalent when they have the same probability function. All 24 optimal sets are equivalent (probability function $pq^3+2p^2q^2+p^3q$). Two adequate sets are isomorphic when one set can be obtained from the other set by renumbering the players.
We notice that equivalency  doesn't imply isomorphic behavior.
So in the asymmetric case we have two different optimal solutions. Both solutions have probability $1-pq$ for success. \\

The last point is to convince ourselves  that any adequate set with
\textit{das}$>$4 doesn't yield better solutions. This can be done by running the
program ASG   with $n=4$,  \textit{das}=5,6,\ldots{}16: all adequate sets have probabilities greater than $pq^3+2{\ p}^2q^2+{
p}^3q$. 
\end{proof}

\section{Symmetric  four person two color Hat Game}
\begin{theorem}
For symmetric four person two color Hat Game we have: maximal probability is $\frac{3}{4}$ with 5 optimal (non isomorphic) decision matrices:
\begin{center}
     \begin{tabular}{|cccccccc|}
     \hline
   {000} & 001     & 010    & 011    & 100   & 101   & 110   & {111} \\
\hline
          & 1     &       & 1     & 0     &       & 0     &     \\
          & 1     &       & 1     & 0     &       & 0     &     \\
          & $\star$      &       &       &       &       &$\star$     &     \\
    0     & 0     &       &       &       &       & 1     & 1 \\
 \hline
    \end{tabular}%
 \end{center}
  \begin{center}
     \begin{tabular}{|cccccccc|}
     \hline
   {000} & 001     & 010    & 011    & 100   & 101   & 110   & {111} \\
\hline
          & 1     & 0      &     &      &    0   & 1     &     \\
          & 1     & 0      &      &      & 0      & 1     &     \\
          & 1      &  0     &       &   0    &       &     &  1   \\
    0     &      &       &   1    &       & 1      & 0     & \\
 \hline
    \end{tabular}%
 \end{center}
 \begin{center}
     \begin{tabular}{|cccccccc|}
     \hline
   {000} & 001     & 010    & 011    & 100   & 101   & 110   & {111} \\
\hline
          & 1     &   1    &    & 0     &       &      &    0 \\
          & 1     &   1    &     & 0     &       &      &    0 \\
        0  & 1     &       &       &       &       &1 &   0  \\
    0     & 1    &       &       &       &       & 1     & 0 \\
 \hline
    \end{tabular}%
 \end{center}

 \begin{center}
     \begin{tabular}{|cccccccc|}
     \hline
   {000} & 001     & 010    & 011    & 100   & 101   & 110   & {111} \\
\hline
    1      &      &       & 1     &     &  0     & 0     &     \\
      1    &      &       & 1     &      &  0     &     0 &     \\
      1    & 0      &       &       &       &       &0     &   1  \\
    1     & 0     &       &       &       &       & 0     & 1 \\
 \hline
    \end{tabular}%
 \end{center}
 
 \begin{center}
     \begin{tabular}{|cccccccc|}
     \hline
   {000} & 001     & 010    & 011    & 100   & 101   & 110   & {111} \\
\hline
      1    & 1     &       &      &      &       & 0     &    0 \\
      1    & 1     &       &      &      &       & 0     &    0 \\
  1        & 1     &       &       &       &       &  0   &    0 \\
    $\star$     &      &       &       &       &       &      & $\star$ \\
 \hline
    \end{tabular}%
 \end{center}
 (Remark: When taking 'pass' for $\star$ in the last matrix, we get the solution where one player passes and the other three go for the well known solution of the three person game)
\end{theorem}
\begin{proof}
Appendix \ref{appendix:three} gives an overview of 40 optimal adequate sets. There are 5 non-isomorphic sets. The first two base decision matrices (corresponding to row 1 and row 13) are  given in section \ref{threetwo}. The last three decision matrices, corresponding to rows 25, 31 and 37 are generated by
adequate sets \{1 2 12 15\}
,
  \{0 3 13 14\} and
  \{0 1 14 15\}.
 They can be found with DMG .
\end{proof}

\section{Comparison between three and four person Hat Game}
In the classic (symmetric) game the three and four player games have the same optimal probability $\frac{3}{4}$: adding a player in the three person game doesn't yield a better result in the classic game. A rather surprising result. In the general version the situation is as follows:
\begin{theorem} The four player game is superior to the three player game when\\ $ \frac{1}{q_2}+\frac{1}{q_3}-1 < \frac{1}{q_1} < \frac{1}{q_2}+\frac{1}{q_3}+\frac{1}{q_4}-1 $ \\ [1 ex]
Maximal difference between optimal probabilities for three and four person game is $q_1q_2q_3p_4$.
\end{theorem}
\begin{proof}
Three person game:\\
CASE: 
$ \frac{1}{q_1} \leq \frac{1}{q_2}+ \frac{1}{q_3}$:\\
optimal probability is 
$p_1+q_1q_2q_3(\frac{1}{q_2}+ \frac{1}{q_3}-\frac{1}{q_1})$.\\
CASE: 
$ \frac{1}{q_1} \geq \frac{1}{q_2}+ \frac{1}{q_3}$:\\
optimal probability is $p_1$ and players 2 and 3 can always pass.\\
\\
Four person game:\\
CASE: $\frac{1}{q_1} \leq \frac{1}{q_2}+\frac{1}{q_3}-1$ \\ maximal winning probability is $p_1+q_1q_2q_3(\frac{1}{q_2}+ \frac{1}{q_3}-\frac{1}{q_1})$ and player 4 can always pass.\\
CASE $ \frac{1}{q_2}+\frac{1}{q_3}-1 \leq \frac{1}{q_1} \leq \frac{1}{q_2}+\frac{1}{q_3}+\frac{1}{q_4}-1 $ \\ maximal winning probability is $p_1+q_1q_2q_3q_4(\frac{1}{q_2}+ \frac{1}{q_3}+\frac{1}{q_4}-\frac{1}{q_1}-1)$\\
CASE $  \frac{1}{q_1} \geq \frac{1}{q_2}+\frac{1}{q_3}+\frac{1}{q_4}-1 $ \\ maximal winning probability is $p_1$ and players 2, 3 and 4 can always pass.\\ \\
The difference  in maximal winning probability between three and four person hat games is:\\
$q_1q_2q_3(1-q_4)(1+\frac{1}{q_1}-\frac{1}{q_2}- \frac{1}{q_3})$ when $\frac{1}{q_2}+ \frac{1}{q_3}-1 < \frac{1}{q_1} \leq \frac{1}{q_2}+ \frac{1}{q_3}$\\
and\\
$q_1q_2q_3q_4(\frac{1}{q_2}+ \frac{1}{q_3}+\frac{1}{q_4}-\frac{1}{q_1}-1)$ when $\frac{1}{q_2}+ \frac{1}{q_3} \leq \frac{1}{q_1} < \frac{1}{q_2}+ \frac{1}{q_3}+\frac{1}{q_4}-1$ \\
Both differences are between 0 and $q_1q_2q_3p_4$
\end{proof}

\section{Computational complexity}
{\raggedright
We consider the number of strategies to be examined to solve the hat problem
with $N$ players and two colors. Each of the $N$ players has $2^{N-1}$ possible
situations to observe and in each situation there are three possible guesses:
white, black or pass. So we have ${(3^{2^{N-1}})}^N$ possible strategies.
Krzywkowski [14] shows that is suffices to examine ${(3^{2^{N-1}-2})}^N$
strategies.
}

{\raggedright
The adequate set method has to deal where \{$i_1,i_2,..,i_{\textit{das}}\}$ with
$0\leq{}i_1<i_2<..<i_{\textit{das}}\leq{}2^N-1$.
}

{\raggedright
The number of strategies for fixed \textit{das} is the number of subsets of dimension \textit{das} of
\{0,1,\ldots{},$\ 2^N$ -1\}: $\left(\begin{array}{l}2^N \\
\textit{das}\end{array}\right)$.
}
But we have to test all possible values of \emph{\textit{das}}. So the correct expression is: $\sum_{\emph{\textit{das}}}
\left(\begin{array}{l}2^N \\
\textit{das}\end{array}\right)=2^{(2^N)}$.
{\raggedright
To get an idea of the power of the adequate set method, we compare the number of
strategies (brute force, Krzywkowski and adequate set method):
}
\begin{center}
\renewcommand{\arraystretch}{1.5}
    \begin{tabular}{|c|c|c|c|}
    \hline
    
   {N}  &   ${(3^{2^{N-1}})}^N$    &   ${(3^{2^{N-1}-2})}^N$   &$2^{(2^N)}$  \\
     
   \hline
      4         & 1.80E+15 & 2.80E+11 &65536 \\
        \hline
    \end{tabular}%
 \end{center}
 
\appendix
        \section{ASG three persons}
\label{appendix:one}      
\includegraphics [scale=0.75]{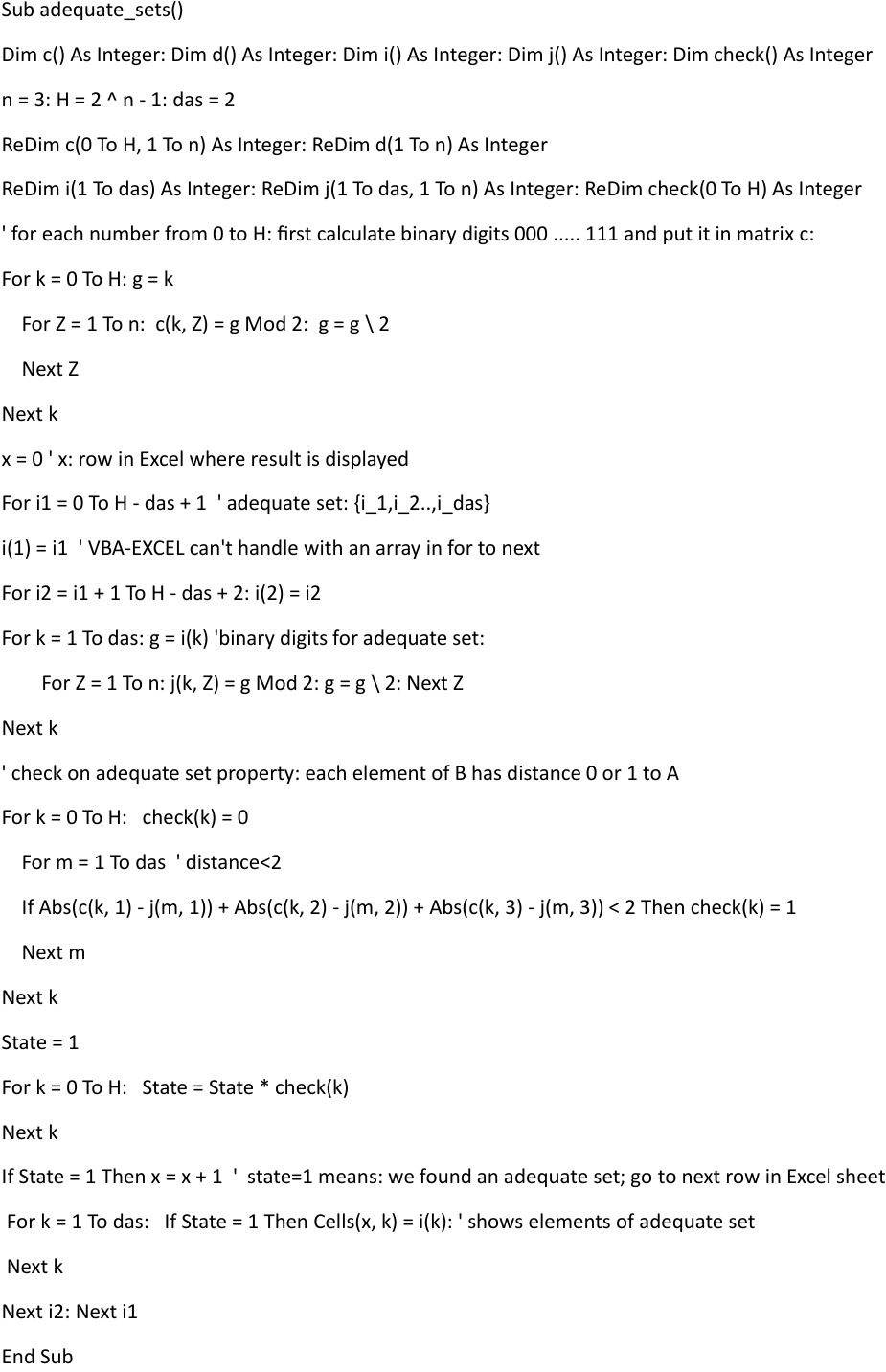}
       \section{DMG three persons}
\label{appendix:two}      
\includegraphics [scale=0.8]{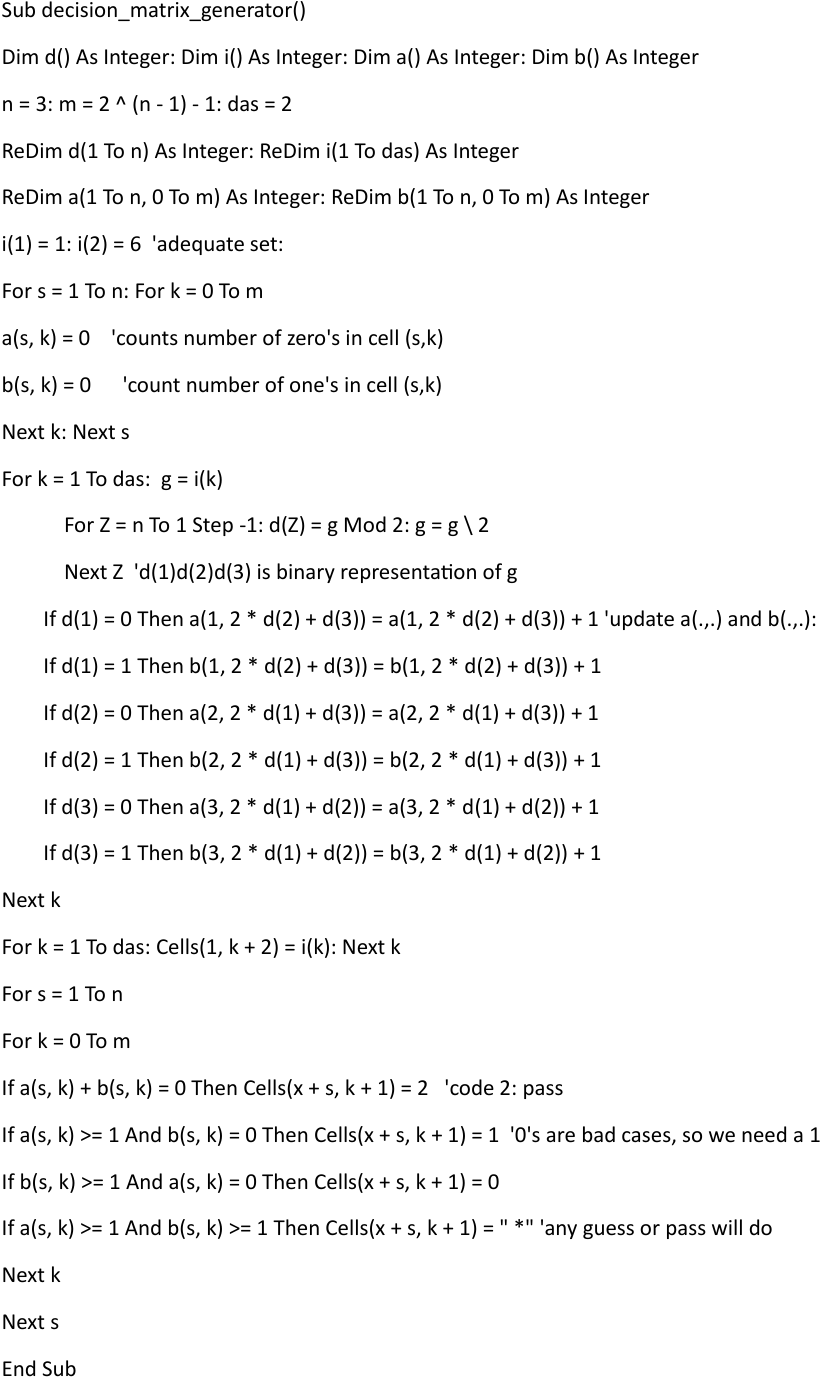}
 \section{Adequate sets; three players; das=3}
 \label{appdas3}
\includegraphics {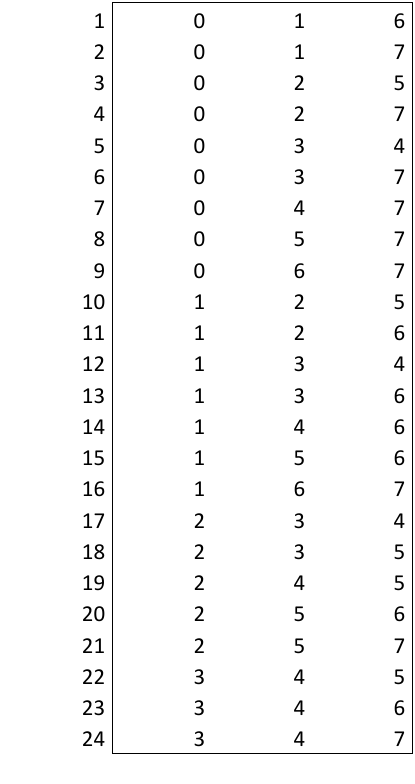}

\section{Dominated adequate sets; General case (three persons, two colors, das=4)}
\label{appendix:seven}
\begin{center}
     \begin{tabular}{|rrrr|c|l|rrrr|c|}
    \cline{1-5}
\cline{7-11}

     
    
    \cline{1-5}
\cline{7-11}

        1     & 3     & 5     & 7     & $q_3$ & &
    1     & 3     & 6     & 7     & $>P(\{1,6\})$ \\
    1     & 5     & 6     & 7     & $>P(\{1,6\})$ &&
    2     & 3     & 5     & 7     & $>P(\{2,5\})$ \\
    2     & 3     & 6     & 7     & $q_2$  &&
    2     & 5     & 6     & 7     & $>P(\{2,5\})$ \\
    3     & 4     & 5     & 7     & $>P(\{3,4\})$ &&
    3     & 4     & 6     & 7     & $>P(\{3,4\})$ \\
    4     & 5     & 6     & 7     & $q_1$  &&
    1     & 3     & 5     & 6     & $>P(\{1,6\})$ \\
    2     & 3     & 5     & 6     & $>P(\{2,5\})$ &&
    3     & 4     & 5     & 6     & $>P(\{3,4\})$ \\
    1     & 2     & 5     & 7     & $>P(\{2,5\})$ &&
    1     & 2     & 6     & 7     & $>P(\{1,6\})$ \\
    1     & 3     & 4     & 7     & $>P(\{3,4\})$   &&
    1     & 4     & 6     & 7     & $>P(\{1,6\})$   \\
    2     & 3     & 4     & 7     & $>P(\{3,4\})$   &&
    2     & 4     & 5     & 7     & $>P(\{2,5\})$   \\
    1     & 2     & 3     & 5     & $>P(\{2,5\})$   &&
    1     & 2     & 3     & 6     & $>P(\{1,6\})$   \\
    1     & 2     & 5     & 6     & $>P(\{1,6\})$   &&
    1     & 3     & 4     & 5     & $>P(\{3,4\})$  \\
    1     & 3     & 4     & 6     & $>P(\{3,4\})$   &&
    1     & 4     & 5     & 6     & $>P(\{1,6\})$   \\
    2     & 3     & 4     & 5     & $>P(\{3,4\})$   &&
    2     & 3     & 4     & 6     & $>P(\{3,4\})$   \\
    2     & 4     & 5     & 6     & $>P(\{2,5\})$   &&
    1     & 2     & 4     & 7     & $\geq q_1$  \\
    1     & 2     & 4     & 5     & $>P(\{2,5\})$   &&
    1     & 2     & 4     & 6     & $>P(\{1,6\})$   \\
    1     & 2     & 3     & 4     & $>P(\{3,4\})$   &&
    0     & 3     & 5     & 7     & $>P(\{0,7\})$   \\
    0     & 3     & 6     & 7     & $>P(\{0,7\})$   &&
    0     & 5     & 6     & 7     & $>P(\{0,7\})$   \\
    0     & 3     & 5     & 6     & $\geq q_1$  &&
    0     & 1     & 3     & 7     & $>P(\{0,7\})$   \\
   0     & 1     & 5     & 7     & $>P(\{0,7\})$   &&
  0     & 1     & 6     & 7     & $>P(\{0,7\})$   \\
   0     & 2     & 3     & 7     & $>P(\{0,7\})$   &&
   0     & 2     & 5     & 7     & $>P(\{0,7\})$   \\
   0     & 2     & 6     & 7     & $>P(\{0,7\})$   &&
    0     & 3     & 4     & 7     & $>P(\{0,7\})$   \\
    0     & 4     & 5     & 7     & $>P(\{0,7\})$   &&
    0     & 4     & 6     & 7     & $>P(\{0,7\})$   \\
    0     & 1     & 3     & 6     & $>P(\{1,6\})$   &&
    0     & 1     & 5     & 6     & $>P(\{1,6\})$   \\
    0     & 2     & 3     & 5     & $>P(\{2,5\})$   &&
    0     & 2     & 5     & 6     & $>P(\{2,5\})$   \\
    0     & 3     & 4     & 5     & $>P(\{3,4\})$   &&
    0     & 3     & 4     & 6     & $>P(\{3,4\})$   \\
    0     & 1     & 2     & 7     & $>P(\{0,7\})$   &&
    0     & 1     & 4     & 7     & $>P(\{0,7\})$   \\
    0     & 2     & 4     & 7     & $>P(\{0,7\})$   &&
    0     & 1     & 2     & 3     & $p_1$ \\
    0     & 1     & 2     & 5     & $>P(\{2,5\})$   &&
    0     & 1     & 2     & 6     & $>P(\{1,6\})$   \\
    0     & 1     & 3     & 4     & $>P(\{3,4\})$   &&
    0     & 1     & 4     & 5     & $p_2$  \\
    0     & 1     & 4     & 6     & $>P(\{1,6\})$   &&
    0     & 2     & 3     & 4     & $>P(\{3,4\})$   \\
    0     & 2     & 4     & 5     & $>P(\{2,5\})$   &&
    0     & 2     & 4     & 6     & $p_3$  \\
\cline{1-5}
\cline{7-11}

   \end{tabular}%
   \end{center}

We have 62 adequate sets, where 54 are absolutely dominated by  $ \{0,7\}$, $\{1,6\}$, $\{2,5\}$ or $\{3,4\}$. \\
Probability of  $\{4,5,6,7\}$ is $q_1p_2p_3+q_1p_2q_3+q_1q_2p_3+q_1q_2q_3=q_1$.\\
For the set $\{0, 3, 5, 6\}$, we obtain that its probability is
$p_1p_2p_3 + p_1q_2q_3 + q_1p_2q_3 + q_1q_2p_3
= q_1(p_2q_3 + q_2p_3 + p_2p_3 + q_2q_3) + (p_1-q_1)(p_2p_3 + q_2q_3)
= q_1 + (p_1-q_1)(p_2p_3 + q_2q_3)
\geq q_1$.\\
Similarly, for the set $\{1, 2, 4, 7\}$, its probability verifies
$p_1p_2q_3 + p_1q_2p_3 + q_1p_2p_3 + q_1q_2q_3
= q_1(p_2p_3 + q_2q_3 + p_2q_3 + q_2p_3) + (p_1-q_1)(p_2q_3 + q_2p_3)
= q_1 + (p_1-q_1)(p_2q_3 + q_2p_3)
\geq q_1.$

\section{four persons, das=4, }
\label{6789}
\includegraphics[scale=0.7]{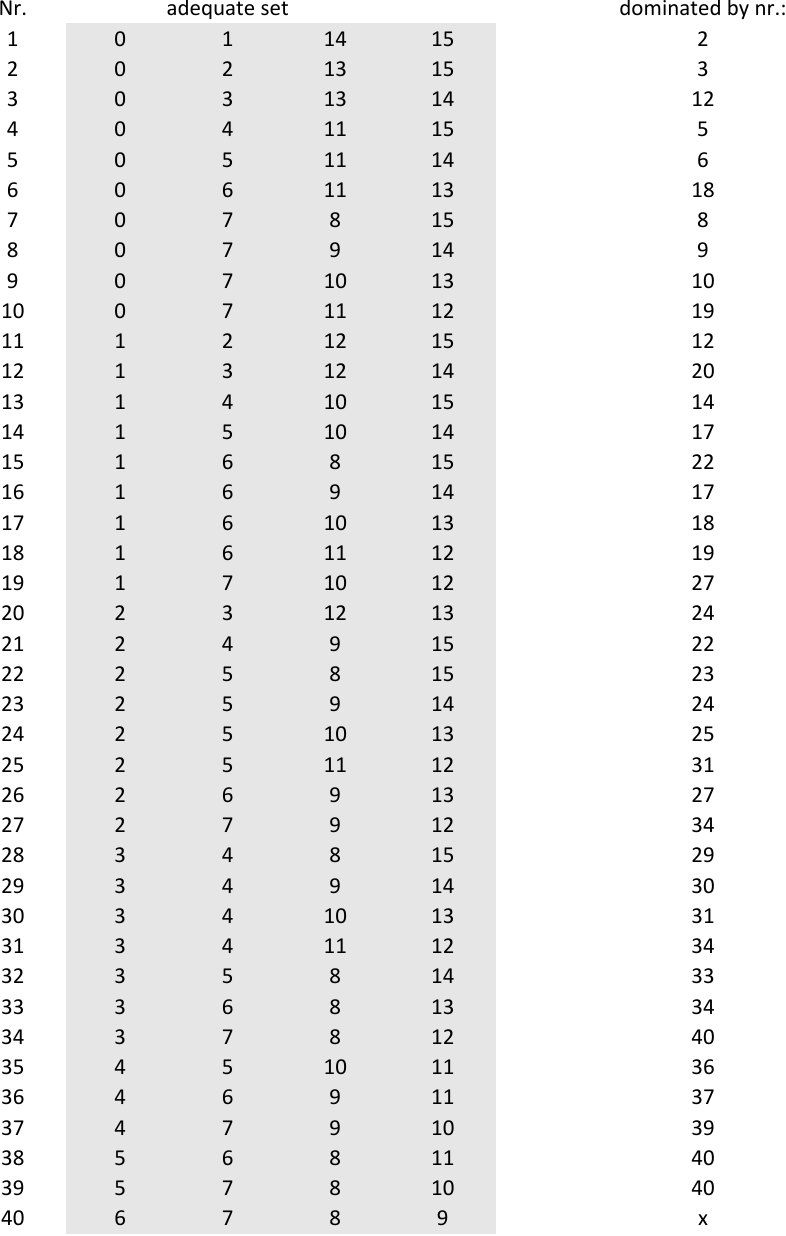}

\includepdf[
    pages=1,
    scale=0.7,
    nup=1x1,
    frame,
    offset={2.5cm, -1.0cm},
    pagecommand={%
        \section{Program: four persons, das=4 versus das=5}
\label{appendix:four}
}       
]{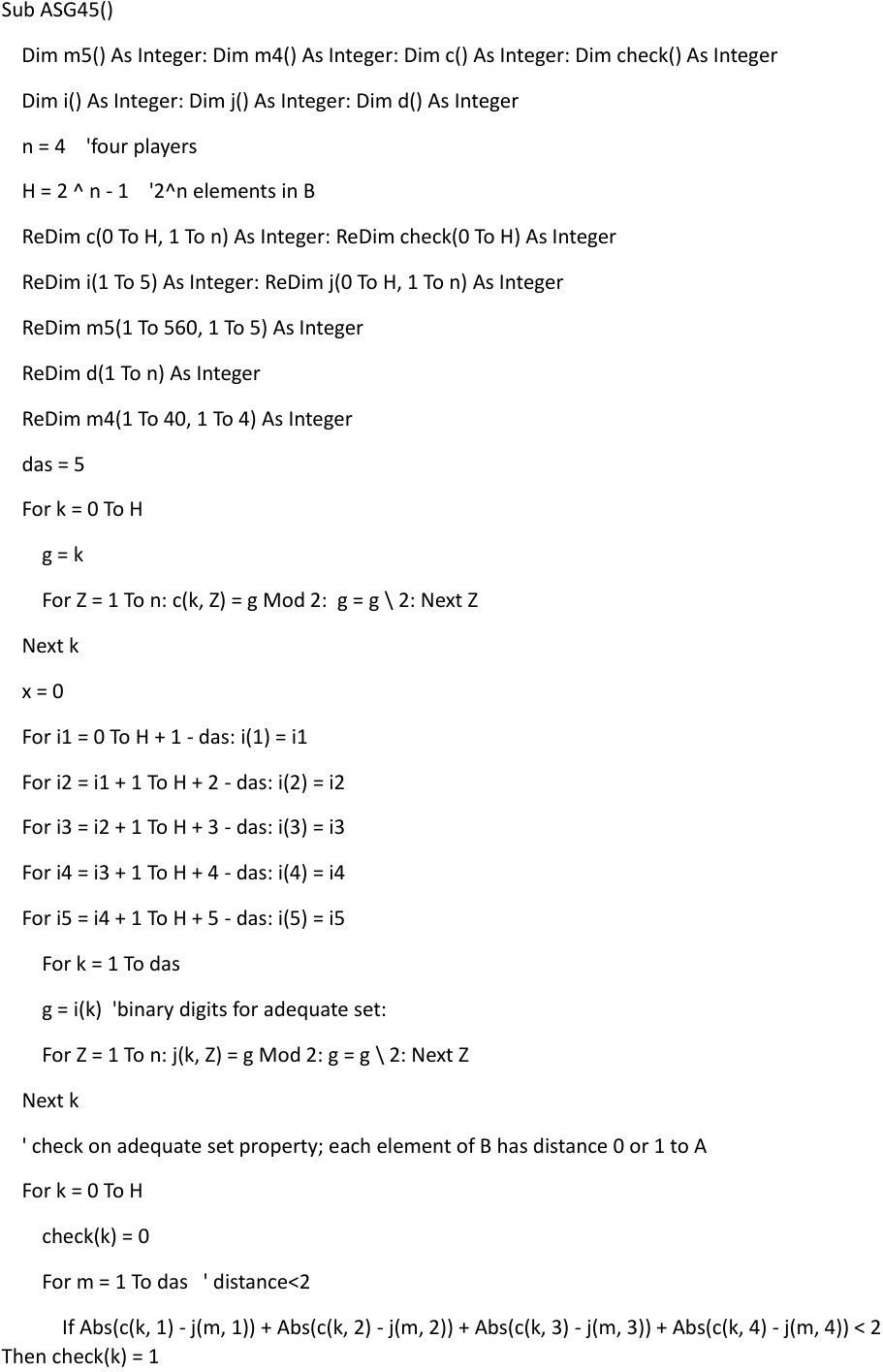}
\includepdf[
    pages=2-3,
    scale=0.7,
    nup=1x1,
    frame,
    offset={2.5cm, -1.0cm},
    pagecommand={%
}       
]{ASG45_program.pdf}

\includepdf[
    pages=1,
    scale=0.7,
    nup=1x1,
    frame,
    offset={2.5cm, -1.0cm},
    pagecommand={%
        \section{Output: four persons, das=4 versus das=5}
\label{appendix:five}
}       
]{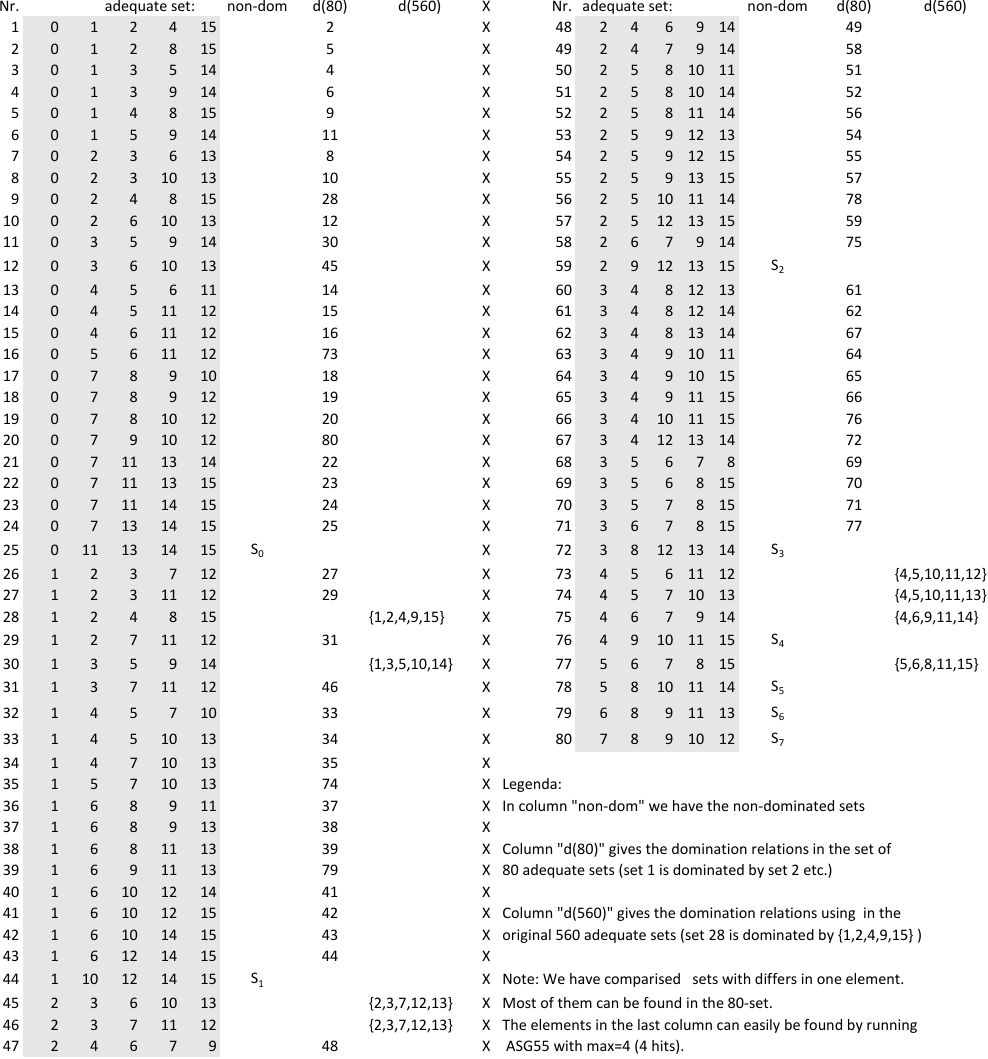}

\section{four persons, das=4, sorted list, p=0.9}
\label{appendix:three}
\includegraphics[scale=0.7]{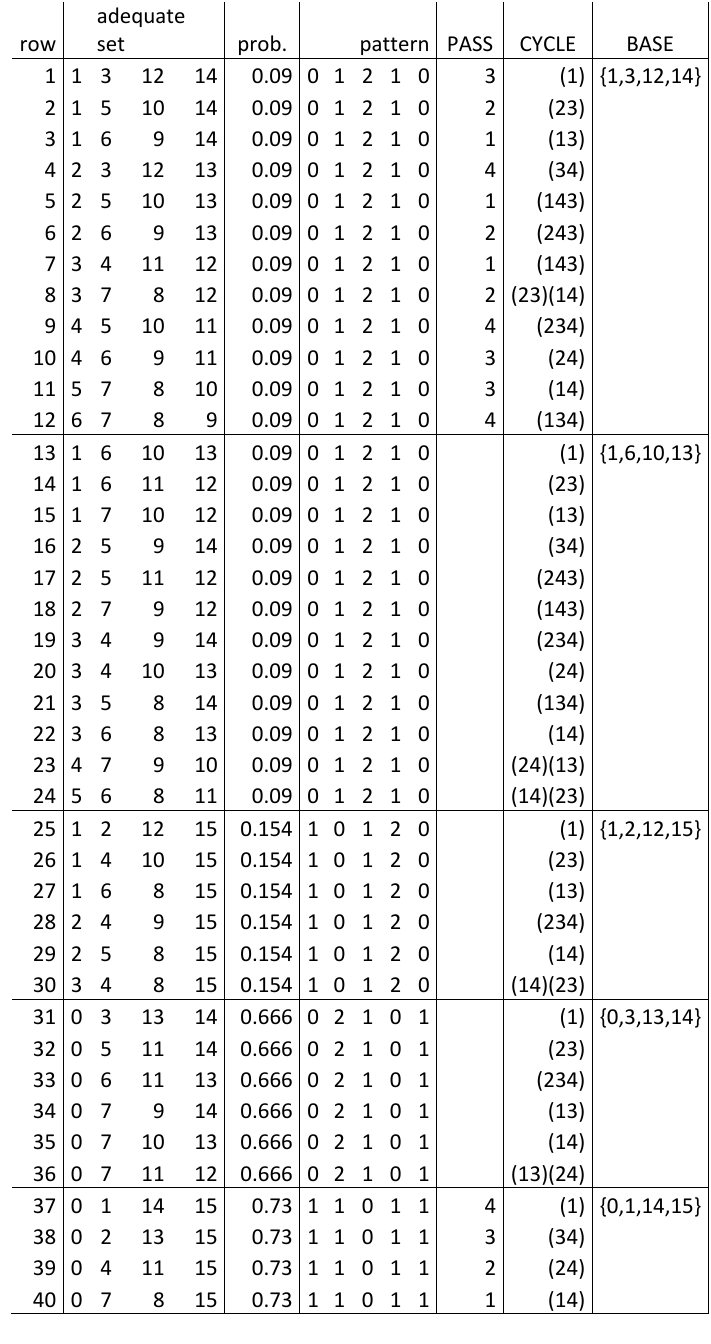}
\clearpage

\end{document}